\documentclass[reqno]{amsart}
\usepackage[ansinew]{inputenc}
\usepackage[all]{xy}
\usepackage{amsmath,latexsym,amssymb,verbatim,amsthm}

\input arrow.tex

\newtheorem{lema}{Lemma}[section]
\newtheorem{teo}[lema]{Theorem}
\newtheorem{pro}[lema]{Proposition}
\newtheorem{defi}[lema]{Definition}
\newtheorem{cor}[lema]{Corollary}
\newtheorem{com}[lema]{Remark}
\newtheorem{eje}[lema]{Example}
\newtheorem{den}[lema]{Notation}
\newtheorem{teo*}{Theorem}

\begin{document}

\title{Reynolds Operator on functors}

\author{Amelia \'{A}lvarez}
\address{Departamento de Matem\'{a}ticas, Universidad de Extremadura,
Avenida de Elvas s/n, 06071 Badajoz, Spain}
\email{aalarma@unex.es}

\author{Carlos Sancho}
\address{Departamento de Matem\'{a}ticas, Universidad de Salamanca, Plaza de la Merced 1-4, 37008 Salamanca, Spain}
\email{mplu@usal.es}

\author{Pedro Sancho}
\address{Departamento de Matem\'{a}ticas, Universidad de Extremadura, Avenida de Elvas s/n, 06071 Badajoz, Spain} \email{sancho@unex.es}



\maketitle

\begin{abstract}
Let $G= {\rm Spec}\, A$ be an affine $R$-monoid scheme. We prove that the category of dual functors (over the category of commutative $R$-algebras) of $G$-modules is equivalent to the category of dual functors of ${\mathcal A}^*$-modules. We prove that $G$ is invariant exact if and only if $A^*= R \times B^*$ as $R$-algebras and the first projection $A^* \to R$ is the unit of $A$. If $\mathbb M$ is a dual functor of $G$-modules and $w_G := (1,0) \in R \times B^* = A^*$, we prove that $\mathbb M^G = w_G \cdot \mathbb M$ and $\mathbb M = w_G \cdot \mathbb M \oplus (1-w_G) \cdot \mathbb M$; hence, the Reynolds operator can be defined on $\mathbb M$.
\end{abstract}

\section*{Introduction}

Let $R$ be a commutative ring with unit. An $R$-module $M$ can be considered as a functor of $\mathcal R$-modules over the category of commutative $R$-algebras, which we will denote by ${\mathcal M}$ and we will call a quasi-coherent module, by defining ${\mathcal M}(S):= M \otimes_R S$. If $\mathbb M$ and $\mathbb N$ are functors of $\mathcal R$-modules, we will denote by ${\mathbb Hom}_{\mathcal R}(\mathbb M,\mathbb N)$ the functor of $\mathcal R$-modules $${\mathbb Hom}_{\mathcal R}(\mathbb M,\mathbb N)(S):={\rm Hom}_{\mathcal S}(\mathbb M_{|S}, \mathbb N_{|S})$$ where $\mathbb M_{|S}$ is the functor $\mathbb M$ restricted to the category of commutative $S$-algebras. The functor $\mathbb M^*:={\mathbb Hom}_{\mathcal R}(\mathbb M,\mathcal R)$ is said to be a dual functor.  For example, ${\mathcal M}$, ${\mathcal M}^*$ and ${\mathbb Hom}_{\mathcal R} ({\mathcal M}, {\mathcal M}')$ are dual functors (see \cite[1.10]{Amel} and Proposition \ref{2.2}).

An affine $R$-monoid scheme $G= {\rm Spec}\,A$ can be considered as a functor of monoids over the category of commutative $R$-algebras: $G^\cdot(S) := {\rm Hom}_{R-alg}(A, S)$, which is known as the functor of points of $G$.

It is well known that the theory of linear representations of algebraic groups can be developed, via their associated functors, as a theory  of (abstract) groups and their linear representations. That is, the category of modules is equivalent to the category of quasi-coherent modules and the category of rational $G$-modules is equivalent to the category of quasi-coherent $G^\cdot$-modules.
Moreover, it is sometimes necessary to consider some natural vector spaces via their associated functors. For example, if $M$ and $M'$ are two (rational) linear representations of $G={\rm Spec}\, A$, then ${\rm Hom}_R(M,M')$ is not a (rational) linear representation of $G$, although $G^\cdot $ operates naturally on ${\mathbb Hom}_{\mathcal R}({\mathcal M}, {\mathcal M}')$; $A^*$ is not a (rational) linear representation of $G={\rm Spec}\, A$, although $G^\cdot $ operates naturally on $\mathcal A^*$.

${\mathcal A}^*$ is a functor of algebras and there exists a natural and obvious morphism of functors of monoids $G^\cdot \hookrightarrow \mathcal A^*$. Then every functor of $\mathcal A^*$-modules is a functor of $G^\cdot$-modules (see Definition \ref{2.1}). In this paper, we prove the following theorem.

\begin{teo*}
The category of dual functors of $G^\cdot$-modules is equal  to the
category of dual functors of ${\mathcal A}^*$-modules.
\end{teo*}

``In the classical situation" observe that if we consider the rational points of $G$, $G^\cdot(R)$, and the inclusion
$G^\cdot(R)\hookrightarrow A^*$, the category of (rational) $G^\cdot(R)$-modules it is not equivalent to the category of $A^*$-modules, even if $G$ is a smooth algebraic group over an algebraically closed field (in this last case it is necessary to introduce topologies on $A^*$ and on the modules).

We say that $G={\rm Spec }\, A$ is invariant exact if taking invariants is an exact functor (see \ref{2.22} and \ref{2.13}). If $A$ is a projective $R$-module and $G$ is invariant exact then it is linearly reductive (\ref{3.5})

We prove that an affine $R$-group scheme $G = {\rm Spec}\, A$ is invariant exact if and only if ${\mathcal A}^* = {\mathcal R} \times {\mathcal B}^*$ as functors of $\mathcal R$-algebras (Corollary \ref{2.9}). If $A$ is a projective $R$-module, $G = {\rm Spec}\, A$ is invariant exact if and only if $A^* = R \times C$ as $R$-algebras (Remark \ref{remark}). If $G$ is invariant exact there exists an isomorphism $A^* = R \times B^*$ such that the first projection $A^* \to R$ is the unit of $A$. Let $w_G:=(1,0) \in {R} \times {B}^*=A^*$ be the ``\textit{invariant integral}'' on $G$. We prove the following theorem.

\begin{teo*} \label{02}
Let $G={\rm Spec}\, A$ be an invariant exact $R$-group and let $w_G \in A^*$ be the invariant integral on $G$. Let $\mathbb M$ be a dual functor of ${G^\cdot}$-modules. It holds that:
\begin{enumerate}
\item $\mathbb M^{G^\cdot} = w_G \cdot \mathbb M$.

\item $\mathbb M$ splits uniquely as a direct sum of $\mathbb M^{G^\cdot}$ and another subfunctor of ${G^\cdot}$-modules, explicitly $$\mathbb M = w_G \cdot \mathbb M \oplus (1-w_G) \cdot \mathbb M.$$
\end{enumerate}
\end{teo*}

We call the projection $\mathbb M \to \mathbb M^{G^\cdot} = w_G \cdot \mathbb M$ the Reynolds operator. Taking sections we have $\mathbb M(R)^{G^\cdot}=w_G\cdot \mathbb M(R)$ and the morphism $\mathbb M(R) \to \mathbb M(R)^{G^\cdot}$, $m\mapsto w_G\cdot m$. The classical Reynolds operator is a particular case of Theorem \ref{02}, for $\mathbb M=\mathcal M$.
The previous theorem still holds for any separated functor of ${\mathcal A}^*$-modules (see Definition \ref{3.1}). More generally, for every functor $\mathbb N$ of $G$-modules, we prove that there exists the maximal separated $G$-invariant quotient of $\mathbb N$ and that the dual of this quotient is $\mathbb N^{*{G^\cdot}}$ (Theorem \ref{3.10}).

In \cite{Huah} it is proved that a Reynolds operator can be defined on ${\rm Hom}_B(M,M')$, where $B$ is a $G$-algebra and  $M$ and $M'$ are two ${ B}G$-modules. Obviously $G^\cdot$ operates on ${\mathbb Hom}_{\mathcal B}(\mathcal M,\mathcal M')$ and it is a separated functor by Proposition \ref{3.13}. Hence the Reynolds operator can be defined on ${\rm Hom}_B(M,M')$ by Theorem \ref{02}. This is an example that shows that functorial treatment can clarify some problems.

Let $\chi\colon G\to G_m$ be a multiplicative character. Given a functor of $G$-modules let $\mathbb M^{\chi}$ be the subfunctor of the  $\chi$-semi-invariant elements of $\mathbb M$ (see Definition \ref{d5}). In Section \ref{section5}, we extend the previous theorems about the invariant integral and the Reynolds operator to the semi-invariant integral and the Reynolds $\chi$-operator. We apply these results to prove some results of \cite{Ferrer} about generalized Cayley's $\Omega$-processes in Example \ref{3.7}.

Finally, these results about the invariant integral, the Reynolds operator, etc., can be extended to functors of monoids with a reflexive functor of functions. We will explain it in detail in a next paper.

\section{Preliminary results}

\cite{Amel} is the basic reference for reading this paper.

Let $R$ be a commutative ring (associative with unit). All functors considered in this paper are covariant functors over the category of commutative $R$-algebras (associative with unit). We will say that $\mathbb X$ is a functor of sets (resp. monoids, etc.) if $\mathbb X$ is a functor from the category of  commutative $R$-algebras to the category of sets (resp. monoids and so forth).

Let $\mathcal R$ be the functor of rings defined by ${\mathcal R}(S):=S$ for every commutative $R$-algebra $S$. We will say that a functor of commutative groups $\mathbb M$ is a functor of $\mathcal R$-modules if we have a morphism of functors of sets ${\mathcal R}\times \mathbb M\to \mathbb M$, so that $\mathbb M(S)$ is an $S$-module for every commutative $R$-algebra $S$. We will say that a functor of $\mathcal R$-modules $\mathbb A$ is a functor of $\mathcal R$-algebras if $\mathbb A(S)$ is a $S$-algebra with unit and $S$ commutes with all the elements of $\mathbb A(S)$. If $\mathbb M$ and $\mathbb N$ are functors of $\mathcal R$-modules, we will denote by ${\mathbb Hom}_{\mathcal R}(\mathbb M,\mathbb N)$ the functor of $\mathcal R$-modules $${\mathbb Hom}_{\mathcal R}(\mathbb M,\mathbb N)(S):={\rm Hom}_{\mathcal S}(\mathbb M_{|S}, \mathbb N_{|S})$$ where $\mathbb M_{|S}$ is the functor $\mathbb M$ restricted to the category of commutative $S$-algebras.

The functor $\mathbb M^*:={\mathbb Hom}_{\mathcal R}(\mathbb M,\mathcal R)$ is said to be a dual functor.

\begin{pro}\label{2.2}
If $\mathbb M^*$ is a dual functor of $\mathcal R$-modules and $\mathbb N$ is a functor of $\mathcal R$-modules, then ${\mathbb Hom}_{\mathcal R} (\mathbb N,\mathbb M^*)$ is a dual functor of $\mathcal R$-modules.
\end{pro}

\begin{proof}
Actually, ${\mathbb Hom}_{\mathcal R} (\mathbb N,\mathbb M^*) = {\mathbb Hom}_{\mathcal R} (\mathbb N \otimes \mathbb M,{\mathcal R})$.
\end{proof}

Given an $R$-module $M$, the functor of $\mathcal R$-modules ${\mathcal M}$  defined by ${\mathcal M}(S) := M \otimes_R S$ is called a quasi-coherent $\mathcal R$-module. The functors $M \rightsquigarrow {\mathcal M}$, ${\mathcal M} \rightsquigarrow {\mathcal M}(R)=M$ establish an equivalence between the category of $\mathcal R$-modules and the category of quasi-coherent $\mathcal R$-modules (\cite[1.12]{Amel}). In particular, ${\rm Hom}_{\mathcal R} ({\mathcal M},{\mathcal M'}) = {\rm Hom}_R (M,M')$. The notion of quasi-coherente $\mathcal R$-module is stable under base change $R\to S$, that is, $\mathcal M_{|S}$ is equal to the quasi-coherent $\mathcal S$-module associated to the $S$-module $M\otimes_R S$.

The functor ${\mathcal M}^* = {\mathbb Hom}_{\mathcal R} ({\mathcal M}, {\mathcal R})$ is  called an $\mathcal R$-module scheme. Specifically, ${\mathcal M}^*(S)={\rm Hom}_S(M\otimes_RS,S)={\rm Hom}_R(M,S)$. It is easy to check that given two functors of $\mathcal R$-modules $\mathbb M$ and $\mathbb M'$, then
$$(\mathbb Hom_{\mathcal R}(\mathbb M_1,\mathbb M_2))_{|S}=
\mathbb Hom_{\mathcal S}(\mathbb M_{|S},\mathbb M'_{|S}).$$
In particular, $(\mathcal M^*)_{|S}$ is an $\mathcal S$-module scheme. An $\mathcal R$-module scheme ${\mathcal M}^*$ is a quasi-coherent $\mathcal R$-module if and only if $M$ is a projective $R$-module of finite type (\cite{Amel2}). A basic result says that quasi-coherent modules and module schemes are reflexive, that is, $${\mathcal M}^{**} = {\mathcal M}$$ (\cite[1.10]{Amel}); thus, the functors $\mathcal M\rightsquigarrow \mathcal M^*$ and $\mathcal M^{*}\rightsquigarrow \mathcal M^{**}=\mathcal M$ establish an equivalence between the category of quasi-coherent modules and the category of module schemes. ${\mathcal M}$ and ${\mathcal M}^*$ are examples of dual functors.

Let $X={\rm Spec}\, A$ be an affine $R$-scheme and let $X^\cdot $ be
the functor of points of $X$, that is, the functor of sets $$X^\cdot (S) ={\rm Hom}_{R-{\rm sch}}({\rm Spec}\, S,X)={\rm Hom}_{R-alg} (A,S)$$  ``points of $X$ with values in $S$''. Given another affine scheme $Y={\rm Spec}\ B$, by Yoneda's lemma $${\rm Hom}_{R-sch}(X,Y)={\rm Hom} (X^\cdot ,Y^\cdot ),$$ and $X^\cdot  \simeq Y^\cdot $ if and only if $X \simeq Y$. We will sometimes denote $X^\cdot  =X$.

Let $R\to S$ be a morphism of rings and let $X_S={\rm Spec}\,
(A\otimes_RS)$, then $(X^\cdot)_{|S}=(X_S)^\cdot$.   Observe that
${\rm Hom}(X^\cdot,\mathcal R)={\rm Hom}(X^\cdot,({\rm Spec}\,
R[X])^\cdot)={\rm Hom}_{R-alg}(R[X],A)$ $=A$, then $\mathbb Hom(X^\cdot,\mathcal R)=\mathcal A$. There is a natural morphism
$X^\cdot  \to {\mathcal A}^*$, because $X^\cdot (S) = {\rm
Hom}_{R-alg} (A,S) \subset {\rm Hom}_R (A,S) = {\mathcal A}^*(S)$.

\begin{pro}
Let $\mathbb M^*$ be a dual functor of $\mathcal R$-modules, any morphism of functors  $X^\cdot \to \mathbb M^*$ factorizes via a unique morphism of functors of $\mathcal R$-modules $\mathcal A^*\to \mathbb M^*$.
\end{pro}

\begin{proof}
It is a consequence of the equalities
$$\aligned {\rm Hom}(X^\cdot,\mathbb M^*) & ={\rm Hom}_{\mathcal R}(\mathbb M,\mathbb Hom(X^\cdot,\mathcal R))={\rm Hom}_{\mathcal R}(\mathbb M,\mathcal A)= {\rm Hom}_{\mathcal R}(\mathbb M\otimes_{\mathcal R}\mathcal A^*,\mathcal R) \\ & ={\rm Hom}_{\mathcal R}(\mathcal A^*,\mathbb M^*).\endaligned$$
\end{proof}

Let $G={\rm Spec}\, A$ be an affine $R$-monoid scheme, that is, $G^\cdot$ is a functor of monoids. ${\mathcal A}^*$ is an $\mathcal R$-algebra scheme, that is, besides from being an $\mathcal R$-module scheme it is a functor of $\mathcal R$-algebras. The natural morphism $G^\cdot  \to {\mathcal A}^*$ is a morphism of functors monoids. By \cite[5.3]{Amel}, given any dual functor of algebras $\mathbb B^*$ (that is, a dual functor of $\mathcal R$-modules which is a functor of $\mathcal R$-algebras), then any morphism of functors of monoids $G^\cdot \to \mathbb B^*$ factorizes via a unique morphism of functors of $\mathcal R$-algebras ${\mathcal A}^*\to\mathbb B^*$.

\begin{defi}\label{2.1}
A functor $\mathbb M$ of (left) $G$-modules is a functor of $\mathcal R$-modules endowed with an action of $G$, i.e., a morphism of functors of monoids $G^\cdot  \to {\mathbb End}_{\mathcal R} (\mathbb M)$. A functor $\mathbb M$ of (left) ${\mathcal A}^*$-modules is a functor of $\mathcal R$-modules endowed with a morphism of functors of $\mathcal R$-algebras ${\mathcal A}^* \to {\mathbb End}_{\mathcal R}(\mathbb M)$.
\end{defi}

The functors $M \rightsquigarrow {\mathcal M}$, ${\mathcal M} \rightsquigarrow {\mathcal M}(R)=M$ establish an equivalence between the category of rational $G$-modules and the category of quasi-coherent $G$-modules. The category of quasi-coherent $G$-modules is equal to the category of quasi-coherent $\mathcal A^*$-modules by \cite[5.5]{Amel}.

\begin{den}
For abbreviation, we sometimes use $g \in G$ or $m \in \mathbb M$ to denote $g \in G^\cdot  (S)$ or $m \in \mathbb M(S)$ respectively. Given $m \in \mathbb M(S)$ and a morphism of $\mathcal R$-algebras $S \to T$, we still denote by $m$ its image by the morphism $\mathbb M(S) \to \mathbb M(T)$.
\end{den}

\begin{defi}
Let $\mathbb M$ be a functor of $G$-modules. We define $$\mathbb M(S)^G:= \{ m \in \mathbb M(S), \mbox{ such that } g \cdot m = m \mbox{ for every } g \in G\}\footnote{More precisely, $g \cdot m=m$ for every $g \in G(T)$ and every morphism of $\mathcal R$-algebras $S \to T$.}$$ and we denote by $\mathbb M^G$ the subfunctor of $\mathcal R$-modules of $\mathbb M$ defined by $\mathbb M^G(S) := \mathbb M(S)^G$. We will say that $m\in \mathbb M$ is left $G$-invariant if $m\in \mathbb M^G$.
\end{defi}

If $\mathbb M$ is a functor of $G$-modules (resp. of right $G$-modules), then $\mathbb M^*$ is a functor of right $G$-modules: $w* g:=w(g\cdot -)$, for every $w \in \mathbb M^*$ and $g\in G$ (resp. of left $G$-modules: $g*w := w(-\cdot g)$). If $\mathbb M_1$ and $\mathbb M_2$ are two functors of $G$-modules and $G$ is an affine $R$-group scheme, then ${\mathbb Hom}_{\mathcal R} (\mathbb M_1,\mathbb M_2)$ is a functor of $G$-modules, with the natural action $g*f := g \cdot f ( g^{-1} \cdot -)$, and it holds that $${\mathbb Hom}_{\mathcal R} (\mathbb M_1,\mathbb M_2)^G = {\mathbb Hom}_G(\mathbb M_1,\mathbb M_2) .$$

Let $\mathbb M$ be a functor of $G$-modules, then $(\mathbb M^G)_{|S} = (\mathbb M_{|S})^{G_S}$.

\section{Invariant Exact Monoids}

From now on, through out this paper $G= {\rm Spec}\, A$ is an affine $R$-monoid scheme.

\begin{teo}\label{funduge}
The category of dual functors of $G$-modules is equivalent to the
category of dual functors of ${\mathcal A}^*$-modules.
\end{teo}

\begin{proof}
Let $\mathbb M$ be a dual functor of $\mathcal R$-modules. By  Proposition \ref{2.2}, ${\mathbb End}_{\mathcal R} (\mathbb M)$ is a
dual functor of $\mathcal R$-algebras. Hence ${\rm Hom}_{mon}
(G^\cdot , {\mathbb End}_{\mathcal R} (\mathbb M)) = {\rm Hom}_{\mathcal R-alg} ({\mathcal A}^*, {\mathbb End}_{\mathcal R} (\mathbb M))$, and giving a structure of functor of $G$-modules on
$\mathbb M$ is equivalent to giving a structure of functor of ${\mathcal A}^*$-modules on $\mathbb M$.

Given two dual functors of $G$-modules $\mathbb M$ and $\mathbb M'$, ${\rm Hom}_G (\mathbb M, \mathbb M') = {\rm Hom}_{{\mathcal A}^*} (\mathbb M, \mathbb M')$: observe that given  a morphism of functors of $\mathcal R$-modules $L: \mathbb M \to \mathbb M'$ and $m \in \mathbb M$, the morphism $L_1: G^\cdot  \to {\mathbb M}'$, $L_1(g) := L(gm) - gL(m)$ is null if and only if the morphism $L_2: {\mathcal A}^* \to \mathbb M'$, $L_2(a) := L(am) - a L(m)$ is null.
\end{proof}

\begin{defi} \label{2.22}
An affine $R$-monoid scheme $G={\rm Spec}\, A$ is said to be left invariant exact if for any exact sequence (in the category of functors of $\mathcal R$-modules) of dual functors of left $G$-modules $$0\to \mathbb M_1 \to \mathbb M_2 \to \mathbb M_3 \to 0$$ the sequence $$0\to \mathbb M_1^G \to \mathbb M_2^G \to \mathbb M_3^G \to 0$$ is exact. $G$ is said to be invariant exact if it is left and right invariant exact.
\end{defi}

If $G$ is an affine  $R$-group scheme and it is left invariant exact, then it is right invariant exact since every functor of right $G$-modules $\mathbb M$ can be regarded as a functor of left $G$-modules: $g \cdot m := m \cdot g^{-1}$.

Let $\Theta: G \to {\mathcal R}$, $g \mapsto 1$ be the trivial character, which induces the trivial representation $\Theta : {\mathcal A}^* \to {\mathcal R}$. Observe that $\Theta = 1 \in A$.

\begin{teo}\label{semisimple}
An affine $R$-monoid scheme $G={\rm Spec}\, A$ is invariant exact if and only if ${\mathcal A}^* = {\mathcal R} \times {\mathcal B}^*$ as $\mathcal R$-algebra schemes, where the projection ${\mathcal A}^* \to {\mathcal R}$ is $\Theta$.
\end{teo}

\begin{proof}
Let us assume that $G$ is invariant exact. The projection $\Theta : {\mathcal A}^* \to {\mathcal R}$ is a morphism of left and right $G$-modules (or ${\mathcal A}^*$-modules). Taking left invariants one obtains an epimorphism $\Theta : {\mathcal A}^{*G} \to {\mathcal R}$. Let $w_l\in A^*$ be left $G$-invariant such that $\Theta(w_l)=1$. Likewise, taking right invariants let $w_r\in A^*$ be right $G$-invariant such that $\Theta(w_r)=1$. Then $w=w_l\cdot w_r \in {\mathcal A}^{*}$ is left and right $G$-invariant and $\Theta(w)=1$. Then, $w' \cdot w= w'(1) \cdot w = w \cdot w' $, because $g \cdot w = w = w \cdot g$. Hence, $w$ is idempotent. Therefore, one finds a decomposition as a product of $\mathcal R$-algebra schemes ${\mathcal A}^* = w \cdot {\mathcal A^*} \oplus (1-w) \cdot {\mathcal A}^*$; moreover, the morphism ${\mathcal R} \to {\mathcal A}^*$, $\lambda \mapsto \lambda \cdot w$, is a section of $\Theta$, $w \cdot {\mathcal A}^* = {\mathcal R} \cdot w$ and $\Theta$ vanishes on $(1-w) \cdot {\mathcal A}^*$

Let us assume now that ${\mathcal A}^* = {\mathcal R} \times {\mathcal B}^*$ and $\pi_1 = \Theta$. Let $w=(1,0) \in {\mathcal R} \times {\mathcal B}^* = {\mathcal A}^*$ and let us prove that $G$ is invariant exact.

For any dual functor of $G$-modules $\mathbb M$, let us see that $w \cdot \mathbb M=\mathbb M^G$. One sees that $w \cdot \mathbb M \subseteq \mathbb M^G$, because $g \cdot (w \cdot m) = (g \cdot w) \cdot m = w \cdot m$, for every $g\in G$ and every $m\in \mathbb M$. Conversely, $\mathbb M^G \subseteq w \cdot \mathbb M$: Let
$m \in \mathbb M$ be $G$-invariant. The morphism $G \to \mathbb M$, $g \mapsto g \cdot m = m$, extends to a unique morphism ${\mathcal A}^* \to \mathbb M$. The uniqueness implies that $w' \cdot m = w'(1) \cdot m$ and then $m = w \cdot m \in w \cdot \mathbb M$.

Taking invariants is a left exact functor. If $\mathbb M_2 \to \mathbb M_3$ is a surjective morphism, then the morphism $\mathbb M_2^G \to \mathbb M_3^G$ is surjective because so is the morphism $\mathbb M_2^G = w \cdot \mathbb M_2 \to w \cdot \mathbb M_3 = \mathbb M_3^G$.
\end{proof}

\begin{com}
If a quasi-coherent $\mathcal R$-module ${\mathcal M }$ is isomorphic to a direct product $\mathbb M\times \mathbb N$ of functors of $\mathcal R$-modules, then $\mathbb M$ and $\mathbb N$ are quasi-coherent (specifically, they are the quasi-coherent modules associated to the modules $\mathbb M(R)$ and $\mathbb N(R)$). Dually, if ${\mathcal M }^*$ is isomorphic to a direct
product $\mathbb M\times \mathbb N$ of functors of $\mathcal R$-modules, then $\mathbb M$ and $\mathbb N$ are $\mathcal R$-module schemes. If ${\mathcal A}^* = \mathbb B \times \mathbb C$ as functors of $\mathcal R$-algebras, then $\mathbb B$ and $\mathbb C$ are $\mathcal R$-algebra schemes.
\end{com}

Let $\chi: G\to G_m$ be a multiplicative character and let $\chi: {\mathcal A}^* \to {\mathcal R}$ be the induced morphism of functors of $\mathcal R$-algebras.

\begin{cor}\label{2.8}
An affine $R$-monoid scheme $G={\rm Spec}\, A$ is invariant exact if and only if ${\mathcal A}^* = {\mathcal R} \times {\mathcal B}^*$ as $\mathcal R$-algebra schemes, where the projection ${\mathcal A}^*\to {\mathcal R}$ is $ \chi$.
\end{cor}

\begin{proof}
The character $\chi$ induces the morphism $G\to {\mathcal A}^*$, $g \mapsto \chi(g) \cdot g$, which induces a morphism of $\mathcal R$-algebra schemes $\varphi \colon {\mathcal A}^* \to {\mathcal A}^*$. This last morphism is an isomorphism because its inverse morphism is the morphism induced by $\chi^{-1}$.

The diagram $$\xymatrix{ {\mathcal A}^* \ar[r]^-{\varphi} \ar[rd]^-\chi & {\mathcal A}^* \ar[d]^-\Theta\\ & {\mathcal R}}$$ is commutative. Hence, via $\varphi$, ``${\mathcal A}^* = {\mathcal R} \times {\mathcal B}^*$ as $\mathcal R$-algebra schemes, where the projection ${\mathcal A}^* \to {\mathcal R}$ is $ \Theta$'' if and only if ``${\mathcal A}^* = {\mathcal R} \times {\mathcal B'}^*$ as $\mathcal R$-algebra schemes, where the projection ${\mathcal A}^* \to {\mathcal R}$ is $ \chi$''. Then, Theorem \ref{semisimple} proves this corollary.
\end{proof}

\begin{cor}\label{2.9}
An affine $R$-group scheme $G={\rm Spec}\, A$ is invariant exact if and only if ${\mathcal A}^* = {\mathcal R} \times {\mathcal B}^*$ as $\mathcal R$-algebra schemes.
\end{cor}

\begin{proof}
Assume that ${\mathcal A}^* = {\mathcal R} \times {\mathcal B}^*$ and let $G \hookrightarrow {\mathcal A}^*$, $g\mapsto g$ be the natural morphism. The composite morphism $$G \hookrightarrow {\mathcal A}^* = {\mathcal R} \times {\mathcal B}^*\overset{\pi_1}{\to} {\mathcal R}$$ is a multiplicative character and $\pi_1$ is the morphism induced by this character. Now it is easy to prove that this corollary is a consequence of Corollary \ref{2.8}.
\end{proof}

If $M$ is an $R$-module and the natural morphism $M \to M^{**}$ is injective, for example if $M$ is a projective module, then $${\rm Hom}_{\mathcal R} ({\mathcal M}^*,{\mathcal M'}^*) = {\rm Hom}_{\mathcal R} ({\mathcal M'} , {\mathcal M}) = {\rm Hom}_{R} (M',M) \subseteq {\rm
Hom}_{ R} ( M^*, M'^*).$$

\begin{com} \label{remark}
Assume that $A$ is a projective $R$-module. If $A^* = C_1 \times C_2$ as $R$-algebras then the morphisms ${\mathcal A}^* \to {\mathcal A}^*$, $w \mapsto (1,0) \cdot w - w \cdot (1,0), (0,1) \cdot w - w \cdot (0,1)$ are null and ${\mathcal A}^*
= {\mathcal A}^* \cdot (1,0) \times {\mathcal A}^* \cdot (0,1)$ as functors of $\mathcal R$-algebras. Then, an affine $R$-group scheme $G={\rm Spec}\, A$ is invariant exact if and only if ${A}^* = {R} \times C$ as $R$-algebras.
\end{com}

\begin{teo}\label{2.10}
An affine $R$-group scheme $G= {\rm Spec}\, A$ is invariant exact if and only if there exists a left $G$-invariant $1$-form $w\in A^*$ such that $w(1)=1$. Moreover, $w$ is unique, it is right $G$-invariant and $*(w)=w$ (where $*$ is the morphism induced on ${\mathcal A}^*$ by the morphism $G \to G$, $g \mapsto g^{-1}$).
\end{teo}

\begin{proof}
If $w_l$ is left invariant and $w_l(1)=1$, then $*w=:w_r$ is right invariant, $w:=w_l\cdot w_r$ is left and right invariant and $w(1)=1$. Now we can proceed as in Theorem \ref{semisimple} in order to prove that $G$ is invariant exact.

Let us only prove the last statement. We follow the notation used in the proof of the last theorem. We know that ${\mathcal A}^{*G} =
(1,0) \cdot {\mathcal A}^* = {\mathcal R} \times 0$, then $(1,0): A \to R$ is the only left $G$-invariant linear map $w: A \to R$ such that $(1,0)(1)=1$. As well, $(1,0)$ is right invariant. Finally,
$*(1,0)$ is left invariant and $(*(1,0))(1) = (1,0)(1) = 1$, then
$*(1,0) = (1,0)$.
\end{proof}

\begin{com}
This result can be found, in \cite{Brion} and \cite{carlos}, when $R$ is a field and $G$ is a linearly reductive algebraic group (that is, every rational $G$-module, $M$, is direct sum of irreducible $G$-modules). If $R$ is an algebraically closed field of characteristic zero, then $G$ is a reductive group if and only if $G$ is linearly reductive, by a theorem of H. Weyl. If $R$ is a field of positive characteristic, then the monoid of matrices $M_n(R)$ is not a linearly reductive monoid (there exists rational representations of $M_n(R)$ no completely reducible), however $M_n(R)$ is invariant exact (observe that given $0\in M_n(R)$ and a $M_n(R)$-module $\mathbb M$ then $0\cdot \mathbb M= \mathbb M^{M_n(R)})$).
\end{com}

In the proof of Theorem \ref{semisimple} we have also proved Theorem \ref{2.11} and Theorem \ref{2.14}.

\begin{teo}\label{2.11}
An affine $R$-monoid scheme $G= {\rm Spec}\, A$ is invariant exact if and only if there exists a left and right $G$-invariant $1$-form $w\in A^*$ such that $w(1)=1$.
\end{teo}

\begin{defi}
Let $G={\rm Spec}\, A$ be an invariant exact affine $R$-monoid scheme. The only 1-form $w_G \in {\mathcal A}^*$ that is left and right $G$-invariant and such that $w_G(1)=1$ is called the invariant integral on $G$ (influenced by the theory of compact Lie groups).
\end{defi}

\begin{teo}\label{2.14}
An affine $R$-group scheme $G= {\rm Spec}\, A$ is invariant exact if and only if for every exact sequence (in the category of functors of $\mathcal R$-modules) of $G$-module schemes
$$0 \to {\mathcal M}_1^* \to {\mathcal M}_2^*
\to {\mathcal M}_3^* \to 0$$
the sequence $$0 \to {\mathcal M}_1^{*G} \to
{\mathcal M}_2^{*G} \to {\mathcal M}_3^{*G} \to 0$$ is exact.
\end{teo}


\begin{teo} \label{2.13}
Let $G={\rm Spec\,} A$ be an affine $R$-monoid scheme. Assume that $A$ is a projective $R$-module. $G$ is invariant exact if and only if the functor ``take invariants'' is (left and right) exact on the category of coherent $G$-modules (or equivalently, the category of rational $G$-modules).
\end{teo}

\begin{proof}
Assume that the functor ``take invariants'' is (left and right) exact on the category of quasi-coherent $G$-modules. $\mathcal A^*$ is an inverse limite of quotients $\mathcal B_i$, which are coherent $\mathcal R$-algebras by \cite[4.12]{Amel}. It can be assumed that the morphism $\Theta\colon \mathcal A^*\to \mathcal R$ factorizes via $\mathcal B_i$ for all $i$. Observe that $\mathcal B_i$ are  (left and right) $\mathcal A^*$-modules, then they are $G$-modules. Now, as in Theorem \ref{2.11}, we can prove that
$\mathcal B_i=\mathcal R\times \mathcal B'_i$ as coherent $\mathcal R$-algebras,(where the projection onto the first factor is $\Theta$).Then, taking inverse limit $\mathcal A^*=\mathcal R\times \mathcal B^*$ and, by Theorem \ref{semisimple}, $G$ is invariant exact.
\end{proof}

\section{Reynolds Operator on separated functors}

Let $\mathbb M$ be a functor of $\mathcal R$-modules and let $\mathbb K$ be the kernel of the natural morphism $\mathbb M \to \mathbb M^{**}$. One has that $\mathbb K(S) = \{ m \in \mathbb M(S) : w(m)=0$, for every $w \in \mathbb M^*(T)$ and every morphism of
$R$-algebras $S \to T \}$. Moreover, $(\mathbb M/\mathbb K)^* = \mathbb M^*$ (then $(\mathbb M/\mathbb K)^{**} = \mathbb M^{**}$) and the morphism $\mathbb M/\mathbb K \to (\mathbb M/\mathbb K)^{**}$ is injective.

\begin{defi}\label{3.1}
We will say that $\mathbb M$ is a separated functor of $\mathcal R$-modules if the morphism $\mathbb M \to \mathbb M^{**}$ is injective, that is, $m \in \mathbb M$ is null if and only if $w(m)=0$ for every $w \in \mathbb M^*$.
\end{defi}

Dual functors of $\mathcal R$-modules are separated: Given $0 \neq m \in \mathbb M=\mathbb N^*$ there exists $n\in \mathbb N$ such that $m(n)\neq 0$; if $\tilde n$ is the image of $n$ by the morphism $\mathbb N\to \mathbb N^{**}=\mathbb M^*$, then $\tilde n(m)=m(n)\neq 0$. Every subfunctor of $\mathcal R$-modules of a separated functor of $\mathcal R$-modules is separated.

\begin{pro} \label{32}
Let $G={\rm Spec}\, A$ be an invariant exact $R$-monoid and let $w_G \in {\mathcal A}^*$ be the invariant integral on $G$. Let $\mathbb M$ be a separated functor of ${\mathcal A}^*$-modules. It holds that:
\begin{enumerate}
\item $\mathbb M^G = w_G \cdot \mathbb M$.

\item $\mathbb M$ splits uniquely as a direct sum of $\mathbb M^G$ and another
subfunctor of $G$-modules, explicitly $$\mathbb M = w_G \cdot
\mathbb M \oplus (1-w_G) \cdot \mathbb M.$$
\end{enumerate}
The morphism $\mathbb M \to \mathbb M^G$, $m \mapsto w_G \cdot m$ will be called the Reynolds operator of $\mathbb M$.
\end{pro}

\begin{proof} $ $
\begin{enumerate}
\item One deduces that $w_G \cdot \mathbb M \subseteq \mathbb M^G$, because $g \cdot (w_G \cdot m) = (g \cdot w_G) \cdot m = w_G \cdot m$ for every $g \in G$ and every $m \in \mathbb M$. Conversely, let us see that $\mathbb M^G \subseteq w_G \cdot \mathbb M$. Let $m \in \mathbb M^G$. The morphism $G \to \mathbb M \hookrightarrow
\mathbb M^{**}$, $g \mapsto g \cdot m = m$, extends to a unique morphism ${\mathcal A}^* \to \mathbb M^{**}$. The uniqueness implies that $w' \cdot m = w'(1) \cdot m$ and then $m = w_G \cdot m \in w_G \cdot \mathbb M$.

\item Since ${\mathcal A}^* = w_G \cdot {\mathcal A^*} \oplus (1-w_G) \cdot {\mathcal A}^*$, then $$\mathbb M = {\mathcal A}^* \otimes_{\mathcal A^*} \mathbb M = w_G \cdot \mathbb M
\oplus (1-w_G) \cdot \mathbb M.$$ Let $\mathbb M = \mathbb M^G \oplus \mathbb N$ be an isomorphism of $G$-modules. The $G$-module structure of $\mathbb N$ extends to an ${\mathcal A}^*$-module structure, because $\mathbb N = \mathbb M/\mathbb M^G$. Moreover, $\mathbb N$ is separated because it is a subfunctor of $\mathcal R$-modules of $\mathbb M$. Now, every morphism of $G$-modules between separated ${\mathcal A}^*$-modules is a morphism of ${\mathcal A}^*$-modules, because the morphism between the
double duals is of ${\mathcal A}^*$-modules by Theorem \ref{funduge}.
Thus, multiplying by $w_G$ one concludes that $w_G \cdot \mathbb M = \mathbb M^G \oplus w_G \cdot \mathbb N$; hence, $w_G \cdot \mathbb N = 0$ and $(1-w_G) \cdot \mathbb M = (1-w_G) \cdot \mathbb N = \mathbb N$.
\end{enumerate}
\end{proof}

\begin{pro}
Let $G={\rm Spec}\, A$ be an invariant exact affine $R$-group scheme and let $\mathbb M$ and $\mathbb N$ be dual functors of $G$-modules. If $\pi : \mathbb M \to \mathbb N$ is an epimorphism of functors of $G$-modules and $s: \mathbb N \to \mathbb M$ is a section of functors of $\mathcal R$-modules of $\pi$, then $w_G \cdot s$ is a section of functors of $G$-modules of $\pi$.
\end{pro}

\begin{proof}
Let us consider the epimorphism of functors of $G$-modules (then of ${\mathcal A}^*$-modules) $$\pi_* : {\mathbb Hom}_{\mathcal R} (\mathbb N,\mathbb M) \to {\mathbb Hom}_{\mathcal R} (\mathbb N,\mathbb N),\, f \mapsto \pi \circ f .$$ Then, $\pi \circ (w_G \cdot s) = \pi_* (w_G \cdot s) = w_G \cdot \pi_*(s) = w_G \cdot
{\rm Id} = {\rm Id}$.
\end{proof}

Likewise, it can be proved that if $\mathbb M$ and $\mathbb N$ are functors of $G$-modules, $\mathbb M$ is a dual functor, $i : \mathbb M \to \mathbb N$ is an injective morphism of $G$-modules and $r$ is a retract of functors $R$-modules of $i$, then $w_G \cdot r$ is a retract of functors of $G$-modules of $i$.

\begin{com} \label{3.5}
We shall say a rational $G$-module $M$ is simple if it does not contain any $G$-submodule $M' \varsubsetneq M$, such that $M'$ is a direct summand of $M$ as an $R$-module (this last condition is equivalent to the morphism of functors of $\mathcal R$-modules ${\mathcal M}^* \to {\mathcal M}'^*$ being surjective, see the previous paragraph to \cite[1.14]{Amel}). If $G$ is an invariant exact affine $R$-group scheme, $M$ is a rational $G$-module and it is a noetherian $R$-module, then it is easy to prove, using the previous proposition, that $M$ is a direct sum of simple $G$-modules.
\end{com}

\begin{eje}
Let us give the proof of the famous \textit{Finiteness Theorem of Hilbert}, \cite{Hilbert}, for its simplicity: ``Let $k$ be a field, let $G$ be a linearly reductive affine $k$-group scheme and let
us consider an operation of $G$ over an algebraic variety  $X={\rm
Spec}\, A$. Then $X/ \sim \, := {\rm Spec}\, A^G$ is an algebraic
variety''.

\begin{proof}
Let  $\xi_1, \ldots, \xi_m$ be a system of generators of the
$k$-algebra $A$. Let $V$ be a finite dimensional $G$-submodule of
$A$ which contains $\xi_1, \ldots, \xi_m$. The natural morphism $S^\cdot V \to A$ is surjective. We have to prove that $A^G$ is an algebra of finite type. As $G$ is invariant exact, it is sufficient to prove that $(S^\cdot V)^G = (k[x_1, \ldots, x_n])^G$ is a $k$-algebra of finite type.

Let $I \subset k[x_1, \ldots, x_n]$ be the ideal generated by $(x_1,
\ldots, x_n)^G$. Let $f_1, \ldots, f_r \in (x_1, \ldots, x_n)^G $ be
a finite system of generators of $I$. We can assume $f_i$ are
homogeneous. Let us prove that $k[x_1, \ldots, x_n]^G = k[f_1,
\ldots, f_r]$. Given a homogeneous $h \in k[x_1, \ldots, x_n]^G$, we
have to prove that $h \in k[f_1, \ldots, f_r]$. We are going to
proceed by induction on the degree of $h$. If ${\rm dg}\, h = 0$
then $h \in k \subseteq k[f_1, \ldots, f_r]$. Let ${\rm dg}\, h = d
> 0$. We can write $h = \overset{r}{\underset{i=1}{\sum}} a_i \cdot
f_i$, where $a_i \in k[x_1, \ldots, x_n]$ are homogeneous of degree
$d - {\rm dg}(f_i)$ (which are less than $d$). Then

$$ h = w_G \cdot h = \overset{r}{\underset{i=1}{\sum}} w_G \cdot (a_i \cdot f_i) \overset*= \overset{r}{\underset{i=1}{\sum}} (w_G \cdot a_i) \cdot f_i $$ (Observe in $\overset*=$ that $g\cdot (a_i \cdot f_i)= (g\cdot a_i)\cdot (g\cdot f_i)=(g\cdot a_i)\cdot f_i$, then $w\cdot (a_i \cdot f_i)=(w\cdot a_i)\cdot f_i$ for all $w\in A^*$). By the induction hypothesis  $w_G \cdot a_i \in k[f_1, \ldots, f_r]$ and therefore $h \in k[f_1, \ldots, f_r]$.
\end{proof}
\end{eje}

Let us see more examples where this theory can be applied.

Let $G= {\rm Spec}\, A$ be an affine $R$-group scheme and let $B$ be an $R$-algebra.

\begin{defi}\label{3.11}
We say that $B$ is a $G$-algebra if $G$ acts on $B$ by endomorphisms of $R$-algebras, that is, there exists a morphism of monoids $G^\cdot  \to {\mathbb End}_{\mathcal R-alg} ({\mathcal B})$.
\end{defi}

We will say that a functor of $\mathcal R$-modules $\mathbb M$ is a functor of $\mathcal B$-modules if there exists a morphism of functors of $\mathcal R$-algebras ${\mathcal B}\to \mathbb End_{\mathcal R}(\mathbb M)$.

\begin{defi}\label{3.12}
Let $B$ be a $G$-algebra and $\mathbb M$ a functor of $\mathcal B$-modules. We say that $\mathbb M$ is a ${\mathcal B}G$-module if it has a $G$-module structure which is compatible with the $\mathcal B$-module structure, that is, $$g (b \cdot m) = g(b) \cdot g(m)$$ for every $g\in G$, $b\in {\mathcal B}$ and $m\in \mathbb M$.
\end{defi}

If $\mathbb M$ and $\mathbb N$ are ${\mathcal B}G$-modules, then it is easy to check that ${\mathbb Hom}_{\mathcal B}({\mathbb M}, {\mathbb N})$ is a subfunctor of $G$-modules of ${\mathbb Hom}_{\mathcal R}({\mathbb M}, {\mathbb N})$ and it coincides with the kernel of the morphism of $G$-modules $$\begin{matrix} {\mathbb Hom}_{\mathcal R} ({\mathbb M}, {\mathbb N}) &
\stackrel{\varphi}{\to} & {\mathbb Hom}_{\mathcal R} ({\mathcal B} \otimes_{\mathcal R} {\mathbb M}, {\mathbb N})
\\ L & \mapsto & L_1 - L_2 \end{matrix} $$ where $L_1(b \otimes m)
:= L(b \cdot m)$ and $L_2(b \otimes m) := b \cdot L(m)$.
Therefore, if $\mathbb N$ is a dual functor as well, then ${\mathbb Hom}_{\mathcal B}({\mathbb M}, {\mathbb N})$ is an ${\mathcal A}^*$-module. Moreover, ${\mathbb Hom}_{\mathcal B} ({\mathbb M}, {\mathbb N})$ is separated because it is an $R$-submodule of ${\mathbb Hom}_{\mathcal R} ({\mathbb M}, {\mathbb N})$, and this latter is separated because it is a dual functor. Hence, if $G = {\rm Spec}\, A$ is an invariant exact $R$-group, ${\mathbb Hom}_{\mathcal B} ({\mathbb M},{\mathbb N})^G = w_G \cdot {\mathbb Hom}_{\mathcal B} ({\mathbb M}, {\mathbb N})$. We have proved the following proposition.

\begin{pro}\label{3.13}
Let $\mathbb N$ be a dual functor of ${\mathcal B}G$-modules and let $\mathbb M$ be a functor
of ${\mathcal B}G$-modules. Then:
\begin{enumerate}
\item ${\mathbb Hom}_{\mathcal B} (\mathbb M, \mathbb N)$ is a separated functor of ${\mathcal A}^*$-modules.

\item If $G={\rm Spec}\, A$ is an invariant exact affine $R$-group scheme, then $${\mathbb Hom}_{\mathcal B} (\mathbb M, \mathbb N)^G = w_G \cdot {\mathbb Hom}_{\mathcal B} (\mathbb M, \mathbb N) $$ and $w_G\cdot\colon {\mathbb Hom}_{\mathcal B} (\mathbb M, \mathbb N)\to {\mathbb Hom}_{\mathcal B} (\mathbb M, \mathbb N)^G$ is the Reynolds operator.
\end{enumerate}
\end{pro}

\section{Reynolds Operator on functors}

Let us generalize the Reynolds operator to all functors of
$G$-modules.

Let us assume that $G={\rm Spec}\, A$ is an invariant exact monoid.

 Given a dual functor of $G$-modules $\mathbb M$, the dual morphism of $\mathbb M^G\hookrightarrow \mathbb M$ is the Reynolds operator of $\mathbb M^*$.

  Let $\mathbb N$ be a separated functor of $G$-modules. Let $\mathbb N_1 = \mathbb N \cap (1-w_G) \cdot \mathbb N^{**}$. Then $\mathbb N_1= \{ n \in \mathbb N : w_G \cdot n = 0 \}$, since $(1-w_G) \cdot \mathbb N^{**} =\{ n' \in \mathbb N^{**} : w_G \cdot n' = 0 \}$. One deduces that $\mathbb N_1^G = \mathbb N_1 \cap \mathbb N^{**G} = 0$ and $(\mathbb N/\mathbb N_1)^G = \mathbb N/\mathbb N_1$, because $\mathbb N/\mathbb N_1$ injects into $\mathbb N^{**}/(1-w_G) \cdot \mathbb N^{**} = \mathbb N^{**G}$. Moreover, $(\mathbb N/\mathbb N_1)^* = \mathbb N^{*G}$: $\mathbb N^{*G} = \mathbb N^* \cdot w_G$ vanishes on $\mathbb N_1$, then $\mathbb N^{*G} \subseteq (\mathbb N/\mathbb N_1)^*$, and $(\mathbb N/\mathbb N_1)^*\subseteq \mathbb N^{*G}$.
Therefore, $\mathbb N^{**}\to(\mathbb N/\mathbb N_1)^{**} $ is the Reynolds operator of $\mathbb N^{**}$.

\begin{teo}\label{3.10}
Let $G = {\rm Spec}\, A$ be an invariant exact $R$-monoid, let $\mathbb N$ be a functor of $G$-modules and let $\mathbb N_1 \subset \mathbb N$ be the subfunctor of $G$-modules defined by $\mathbb N_1 := \{ n \in \mathbb N : w_G \cdot \tilde n = 0 \}$, where $\tilde n$ denotes the image of $n$ by the morphism $\mathbb N \to \mathbb N^{**}$. It holds that:
\begin{enumerate}
\item $\mathbb N/\mathbb N_1$ is the maximal separated $G$-invariant quotient of
$\mathbb N$.

\item The double dual of the morphism $\mathbb N \to \mathbb N/\mathbb N_1$ is the Reynolds operator $\mathbb N^{**}\to \mathbb N^{**G}$ and one has the
commutative diagram \begin{equation}\label{a}\xymatrix{ \mathbb N \ar[r] \ar[d] & \mathbb N/\mathbb N_1 \ar@{^{(}->}[d] \\ \mathbb N^{**} \ar[r] & \mathbb N^{**G}=(\mathbb N/\mathbb N_1)^{**}.} \end{equation}

\item If $\mathbb N$ is a dual functor, then $\mathbb N/\mathbb N_1= \mathbb N^G$ and the morphism $\mathbb N \to \mathbb N/\mathbb N_1$ is the Reynolds operator of $\mathbb N$.
\end{enumerate}
\end{teo}

\begin{proof} $ $
\begin{enumerate}
\item[(2)] $\mathbb N_1$ is the kernel of the composite morphism $\mathbb N \to \mathbb N^{**} \to \mathbb N^{**G}$, then $\mathbb N_1$ contains the kernel $\mathbb K$ of the morphism $\mathbb N \to \mathbb N^{**}$. Let $\mathbb N'=\mathbb N/\mathbb K$. Observe that $\mathbb N'^*= \mathbb N^*$, that $\mathbb N'$ is separated,
    $\mathbb N'_1=\mathbb N_1/\mathbb K$ and  $\mathbb N'/\mathbb N'_1=\mathbb N/\mathbb N_1$. Therefore the diagram (\ref{a})
    is commutative because is commutative for
    $\mathbb N=\mathbb N'$. In particular, $\mathbb N/\mathbb N_1$ is separated.

\item[(1)]We must prove that if $\mathbb P \subseteq \mathbb N$ is a subfunctor of $G$-modules such that $\mathbb N/\mathbb P$ is separated and $G$-invariant, then $\mathbb N_1 \subseteq \mathbb P$, i.e., $\mathbb N/\mathbb P$ is a quotient of $\mathbb N/\mathbb N_1$.

$\mathbb N/(\mathbb P \cap \mathbb N_1)$ is $G$-invariant and separated functor, because the morphism $\mathbb N/(\mathbb P \cap \mathbb N_1) \hookrightarrow \mathbb N/\mathbb N_1 \oplus \mathbb N/\mathbb P$, $\bar{h} \mapsto (\bar{h}, \bar{h})$ is injective. It is enough to prove that $\mathbb N_1 = \mathbb P \cap \mathbb N_1$. Let us denote $\mathbb P' = \mathbb P \cap \mathbb N_1$. From the composition of injections $\mathbb N^{*G} = (\mathbb N/\mathbb N_1)^* \hookrightarrow (\mathbb N/\mathbb P')^* \hookrightarrow \mathbb N^{*G}$ one concludes that $(\mathbb N/\mathbb N_1)^* = (\mathbb N/\mathbb P')^*$. Now, the commutative diagram $$\xymatrix{ \mathbb N/\mathbb P' \ar[r] \ar@{^{(}->}[d] & \mathbb N/\mathbb N_1 \ar@{^{(}->}[d] \\ (\mathbb N/\mathbb P')^{**} \ar@{=}[r]  & (\mathbb N/\mathbb N_1)^{**}}$$ implies that the morphism $\mathbb N/\mathbb P' \to \mathbb N/\mathbb N_1$ is injective; hence, $\mathbb P'=\mathbb N_1$.

\item[(3)] Recall that $\mathbb N = w_G \cdot \mathbb N \oplus (1-w_G) \cdot \mathbb N$, $\mathbb N_1 = (1-w_G) \cdot \mathbb N$ and $\mathbb N/\mathbb N_1 = w_G \cdot \mathbb N = \mathbb N^G$.
\end{enumerate}
\end{proof}

Let us observe that $\mathbb N_1^0 := \{ w \in \mathbb N^* : w(\mathbb N_1) = 0 \} = (\mathbb N/\mathbb N_1)^* = \mathbb N^{*G}$. On the other hand, $(\mathbb N^{*G})^0 := \{ n \in \mathbb N
: \mathbb N^{*G}(n) = 0 \} = \{ n \in \mathbb N : \tilde{n} (\mathbb N^{*G} = \mathbb N^* \cdot w_G ) = 0 \} = \{ n \in \mathbb N : (w_G \cdot \tilde{n}) (\mathbb N^{*}) = 0 \} = \{ n\in \mathbb N : w_G \cdot \tilde{n} = 0 \} = \mathbb N_1$.

\section{Semi-invariants}\label{section5}

Let $\chi\colon G={\rm Spec}\,A\to {\mathcal R}$ be a multiplicative character and let $\chi\colon \mathcal A^*\to \mathcal R$ be the induced morphism.

\begin{defi}\label{d5}
Let $\mathbb M$ be a functor of $G$-modules. An element $m \in \mathbb M$ is said to be (left) $\chi$-semi-invariant if $g \cdot m = \chi(g) \cdot m$ for every $g\in G$.
\end{defi}

\begin{defi}
Let $G={\rm Spec}\, A$ be an affine $R$-monoid scheme and let $\chi \colon G\to {\mathcal R}$ be a multiplicative character. We will call
the 1-form $w_\chi \in { A}^*$ which is left and right $\chi$-semi-invariant and such that $w_\chi(\chi)=1$, if it exists, a $\chi$-semi-invariant integral on $G$.
\end{defi}

If a $\chi$-semi-invariant integral exists, then it is unique: Observe that $w \cdot w_\chi = w(\chi) \cdot w_\chi $, because $g \cdot w_\chi = \chi(g) \cdot w_\chi$ for every $g \in G$, since $w_{\chi}$ is left $\chi$-semi-invariant. Likewise, $ w_\chi \cdot w=w(\chi) \cdot w_\chi$. Given a left $\chi$-semi-invariant $w\in A^*$ such that $w(\chi)=1$, one concludes that $w=w_\chi(\chi) \cdot w = w_\chi \cdot w = w(\chi) \cdot w_\chi = w_{\chi}$.

Given a functor $\mathbb M$ of $G$-modules we will define $\mathbb M^\chi$ to be the functor $\mathbb M^\chi (S):= \{ m \in \mathbb M(S) \colon g \cdot m = \chi(g) \cdot m$, for every $g \in G(T)$, and every $S$-algebra $T\}$.

\begin{pro}
Let $G={\rm Spec}\, A$ be an affine $R$-monoid scheme and let $\chi \colon G\to {\mathcal R}$ be a multiplicative character. $\mathcal A^*=\mathcal R\times\mathcal B^*$ as functors of $\mathcal R$-algebras, where the projection onto the first factor is $\chi$
if and only if there exists the $\chi$-semi-invariant integral on $G$.
\end{pro}

\begin{proof}
If there exists the $\chi$-semi-invariant integral on $G$, $w_\chi$, then ${\mathcal A^*} = w_\chi \cdot {\mathcal A^*} \times (1-w_\chi) \cdot {\mathcal A}^*=\mathcal R\times\mathcal B^*$ as functors of $\mathcal R$-algebras, where the projection onto the first factor is $\chi$. Conversely, if ${\mathcal A}^* = {\mathcal R} \times {\mathcal B}^*$, where the first projection ${\mathcal A}^*\to {\mathcal R}$ is $\chi$, then $w_\chi=(1,0) \in {\mathcal R} \times {\mathcal B}^*$.
\end{proof}

By Corollary \ref{2.8} we obtain the following theorem.

\begin{teo}
Let $G$ be an affine monoid scheme and let $\chi\colon G\to G_m$ be a multiplicative character. $G$ is invariant exact if and only if there exists the $\chi$-semi-invariant integral on $G$.
\end{teo}

Likewise as in Proposition \ref{32}, we obtain the following result.

\begin{pro}
Let $G={\rm Spec}\, A$ be an affine $R$-monoid scheme  and assume there exists the $\chi$-semi-invariant integral on $G$, $w_\chi \in
{\mathcal A}^*$. Let $\mathbb M$ be a separated functor of ${\mathcal A}^*$-modules. It holds that:
\begin{enumerate}
\item $\mathbb M^\chi = w_\chi\cdot \mathbb M$.

\item $\mathbb M$ splits uniquely as a direct sum of $\mathbb M^\chi$ and another subfunctor of $G$-modules, explicitly $$\mathbb M=w_\chi \cdot \mathbb M \oplus (1-w_\chi) \cdot \mathbb M.$$
\end{enumerate}
\end{pro}

We call the morphism $\mathbb M\to \mathbb M^{\chi}$, $m\mapsto w_{\chi}\cdot m$, the Reynolds $\chi$-operator.

\begin{eje}\label{3.7}
Let $G={\rm Spec}\, A$ be an affine $R$-monoid scheme and let $\chi \colon
G \to {\mathcal R}$ be a multiplicative character. An $\Omega$-process associated to $\chi$ (see \cite[3.1]{Ferrer}) is a nonzero linear operator $\Omega \colon A \to A$ such that $$\Omega(a\cdot g) =\chi(g)\cdot (\Omega(a)\cdot g); \,\, \Omega(g\cdot a) =\chi(g)\cdot (g\cdot\Omega(a)) $$ for all $a \in A$ and $g \in G$. The composite morphism $$A\overset{\Omega} {\longrightarrow} A \overset{\chi\cdot}{\longrightarrow} A$$ is a morphism of left and right $G$-modules: $(\chi\cdot\circ \Omega)(g \cdot a)=\chi\cdot \chi(g)\cdot (g\cdot \Omega(a))=(g\cdot \chi)\cdot (g\cdot\Omega(a))=g\cdot (\chi\cdot \Omega (a))=g\cdot ((\chi\cdot \circ \Omega)(a))$, likewise we prove that $\chi\cdot\circ \Omega$ is a morphism of right $G$-modules.
Since $${\rm Hom}_{\text{left-right $G$-modules}}(A,A)={\rm
Hom}_{\text{left-right ${\mathcal A}^*$-modules}}({\mathcal A}^*,{\mathcal A^*}) = Z(A^*)$$ ($Z(A^*)$ is the center of $A^*$), then $\chi\cdot\circ \Omega = z\cdot $ for some $z \in Z(A^*)$. If $G$ is a linearly reductive monoid and $R$ is an algebraically closed field, then $\mathcal A^*$ is a semisimple algebra scheme and $A^*=\prod_{E_i\in I} {\rm End}_R (E_i)$ where $I$ is the set of irreducible representations of $G$ (up to isomorphism), by \cite[6.2, 6.8]{Amel}, hence $\chi\cdot\circ \Omega\in Z(A^*)=\prod_{i\in I}R$ (on the other hand, see \cite[4.4]{Ferrer}).

Assume now that $0 \in G$ (that is an element such that $0 \cdot g = g \cdot 0 = 0$ for all $g \in G$) and that $\Omega(\chi)=1$ (generally $\chi\cdot \Omega (\chi)=z\cdot\chi=\chi(z)\cdot \chi=\chi\cdot\chi(z)$, $\chi(z)\in R$). The projection $w:A \to R$, $a\mapsto a(0)$ is left and right invariant and $w(1)=1$, then $G$ is an invariant exact $R$-monoid and $w=w_G$. The composite morphism $w'= w_G \circ \Omega$ is left and right $\chi$-semi-invariant and $w'(\chi)=1$, then $w'=w_{\chi}$. Given a rational $G$-module $M$, let us calculate the Reynolds $\chi$-operator of $M$, that is, the morphism $M\to M$, $m\mapsto w_\chi\cdot m$ (on the other hand, see \cite[5.1]{Ferrer}). The dual morphism of the multiplication morphism ${\mathcal M}^* \otimes {\mathcal A}^*\to  {\mathcal M}^*$ is the comultiplication morphism $\mu \colon M \to M \otimes A$. If $\mu(m)= \sum_l m_l \otimes a_l$, then $g \cdot m = \sum_l a_l(g) \cdot m_l$, for all $g \in G$. Hence,

$$w_\chi\cdot m = \sum_l a_l(w_\chi)\cdot m_l= \sum_l w_\chi(a_l)\cdot m_l= \sum_l \Omega (a_l)(0)\cdot m_l .$$
\end{eje}


\begin{thebibliography}{99}


\bibitem{Amel} \textsc{\'{A}lvarez, A., Sancho, C., Sancho, P.,}
\textit{\!Algebra schemes and their representations}, J. Algebra
{\bf 296/1} (2006) 110-144.

\bibitem{Amel2} \textsc{\'{A}lvarez, A., Sancho, C., Sancho, P.,}
\textit{\!Characterization of Quasi-Coherent Modules  that
are Module Schemes}, Communications in Algebra (2009),37:5,1619 — 1621.


\bibitem{Brion} \textsc{Brion, M., Schwarz, G.W.,}
\!Th{\'{e}}orie des invariants et g{\'{e}}ometrie des
vari{\'{e}}t{\'{e}}s quotients, Travaux en cours, vol. 61,
Hermann, Paris, 2000.

\bibitem{Huah} \textsc{Chu, H., Hu, S.-J., Kang, M.-C.,}
\textit{A variant of the Reynolds operator}, Proc. Amer. Math.
Soc. {\bf 133/10} (2005) 2865-2871.

\bibitem{Ferrer} \textsc{Ferrer Santos, W., Rittatore, A.,}
\textit{\!Generalizations of Cayley $\Omega$-process}, Proc. Amer. Math. Soc. {\bf 135/4} (2007) 961-968.

\bibitem{Hilbert} \textsc{Hilbert, D.,}
\textit{\!\"{U}ber die Theorie der algebraischen Formen}, Math. Ann. {\bf 36} (1890) 473-534.



\bibitem{carlos} \textsc{Sancho de Salas, C.,} \!Grupos
algebraicos y teor\'{\i}a de invariantes, Sociedad Matem\'{a}tica
Mexicana, M\'{e}xico, 2001.


\end{thebibliography}
\end{document}